\documentclass[11pt]{article}

\usepackage[a4paper,margin=1.05in]{geometry}
\usepackage[utf8]{inputenc}
\usepackage[T1]{fontenc}
\usepackage{lmodern}
\usepackage{microtype}
\usepackage{amsmath,amssymb,amsfonts,amsthm,mathtools}
\usepackage{enumitem}
\usepackage{extarrows}
\usepackage{booktabs}
\usepackage{array}
\usepackage{xcolor}
\usepackage{colortbl}
\usepackage{caption}
\usepackage[colorlinks=true,linkcolor=blue,citecolor=blue]{hyperref}

\AddToHook{cmd/normalsize/after}{%
  \addtolength{\abovedisplayskip}{0.75pt}%
  \addtolength{\belowdisplayskip}{0.75pt}%
  \addtolength{\abovedisplayshortskip}{0.5pt}%
  \addtolength{\belowdisplayshortskip}{0.5pt}%
}
\allowdisplaybreaks
\numberwithin{equation}{section}

\newtheorem{theorem}{Theorem}[section]
\newtheorem{proposition}[theorem]{Proposition}
\newtheorem{lemma}[theorem]{Lemma}

\theoremstyle{remark}
\newtheorem{remark}[theorem]{Remark}
\theoremstyle{plain}

\newcommand{\R}{\mathbb R}
\newcommand{\C}{\mathbb C}

\newcommand{\cE}{\mathcal E}
\newcommand{\cH}{\mathcal H}
\newcommand{\cQ}{\mathcal Q}

\newcommand{\eps}{\varepsilon}
\newcommand{\dd}{\,\mathrm d}
\newcommand{\Cat}{\mathrm{Cat}}

\newcommand{\wh}[1]{\widehat{#1}}

\newcommand{\norm}[1]{\left\lVert #1\right\rVert}
\newcommand{\Sym}{\operatorname{Sym}}

\title{A uniform positive lower bound for the analyticity radius of KdV on the real line}
\author{Ming Wang \and Jie Zeng}
\date{}

\begin{document}
\maketitle

\begin{abstract}
We show that, for the KdV equation on the real line
with initial data in an $L^2$-based analytic class, the
radius of spatial analyticity of the solution has a fixed positive
lower bound for all time.
This improves the previous quantitative lower bounds, which tend to
zero as $|t|\to\infty$.
The proof relies on a new positive measure representation of the Hamiltonian
multipliers to sum the KdV Hamiltonians into an exactly conserved
energy that controls the analytic norm.
\end{abstract}



\section{Introduction}\label{section-introduction}

 In this paper, we consider the Cauchy problem for the Korteweg--de Vries (KdV) equation
\begin{equation}\label{equ-kdv}
 \partial_tq+\partial_x^3q-6q\partial_xq=0,\qquad q(0,x)=q_0(x),\qquad (t,x)\in\R^2.
\end{equation}
The KdV equation arose in the study of long water
waves in a shallow channel \cite{KdV}. It is a classical example of a completely
integrable equation. The global well-posedness of the KdV equation
in Sobolev spaces $H^s(\R)$ is now well understood (see
\cite{CKSTT,KV} and the references therein).

In this work, we are interested in the global well-posedness of the
KdV equation in analytic spaces. Our motivation comes from two analytic
smoothing effects of the KdV equation. Indeed, Kato and Ogawa
\cite{KO} showed that the solution with Dirac initial data is jointly
real analytic in space and time for small nonzero times. Tarama
\cite{Tar} and Rybkin \cite{Ryb} obtained spatial analyticity and
meromorphic continuation, respectively, for $t>0$ under suitable
one-sided exponential decay assumptions.

To state the problem precisely, we first introduce the analytic spaces used in this paper. For every $\sigma\geq0$, let $G^\sigma=G^{\sigma}(\R)$ be the Banach space endowed with the norm
\begin{equation}\label{equ-intro-9-22}
	 \norm f_{G^\sigma}:=\norm{e^{\sigma|D|}f}_2=\norm{e^{\sigma|\xi|}\widehat{f}}_2.
\end{equation}
Here $\wh f(\xi)=\int_\R e^{-ix\xi}f(x)\dd x$ denotes the Fourier transform and $\|\cdot\|_2$ denotes the standard norm in $L^2(\R)$.
By the Paley-Wiener theorem, every $f\in G^\sigma$ can be extended to a holomorphic
functions on the strip $\{z\in \C: |\operatorname{Im}z|<\sigma\}$.

The local well posedness in $G^\sigma$ is easy and follows from the classical $L^2$ theory. In fact, if we define the Bourgain space endowed with norm
$$
\|q\|_{X^{0,b}}=\|(1+|\tau-\xi^3|)^b\wh{q}(\tau,\xi)\|_{L^2(\R^2)}
$$
and its variants in $G^\sigma$
$$
\|q\|_{G^{\sigma,b}}=\|e^{\sigma|\xi|}(1+|\tau-\xi^3|)^b\wh{q}(\tau,\xi)\|_{L^2(\R^2)}
$$
then the classical bilinear estimates in $X^{0,b}$ implies\footnote{This is a direct conseqence of the inequality $e^{\sigma|\xi_1+\xi_2|}\leq e^{\sigma|\xi_1|}e^{\sigma|\xi_2|}$.} the bilinear estimates in $G^{\sigma,b}$
\begin{align}\label{equ-bi-G}
\|\partial_x(q_1q_2)\|_{G^{\sigma,b'-1}}\lesssim_b \|q_1\|_{G^{\sigma,b}}\|q_2\|_{G^{\sigma,b}}	
\end{align} 
for $0<b<b'\leq 3/4$. With \eqref{equ-bi-G} in hand, then the contraction mapping principle gives the local well posedness in $G^\sigma$, see e.g. \cite{HW,SdS}.

 For the global well posedness of \eqref{equ-kdv}, it requires to show the persistence: the solution $q(t)\in G^\sigma$ for all time $t\in \R$ if $q_0\in G^\sigma$. The problem is difficult since there is no known conservation law in $G^\sigma$. Thus one may consider the following more accessible problem:
\begin{equation}\label{equ-intro-1128-1}
 \mbox{If $q_0\in G^{\sigma_0}$, how large $\sigma_q(t)$ can be as $t\to \infty$? }
\end{equation} 
Here and below, for a solution $q$ of \eqref{equ-kdv}, we define
its analytic radius by
$$
	\sigma_q(t):=\sup\{\sigma>0:q(t)\in G^\sigma(\R)\}.
$$

There are several progresses toward the question \eqref{equ-intro-1128-1}.  Kato--Masuda \cite{KM} used Lyapunov functions with a varying
analyticity parameter to obtain a positive strip on every finite
time interval.   In
\cite[Theorem~11 and Remark~12]{BG},
Bona--Gruji\'c established the lower bound $\sigma_q(t)\geq e^{-ct^2}$ for
large $|t|$ through Gevrey estimates.  In \cite{BGK}, Bona--Gruji\'c--Kalisch then showed the improved
algebraic bounds $\sigma_q(t)\geq  ct^{-12}$ for large $t$.  In \cite{SdS}, Selberg--da Silva further proved a lower bound of order $|t|^{-(4/3+\eps)}$
for every $\eps>0$. The main idea in \cite{SdS} is as follows. Applying $e^{\sigma|D|}$
to \eqref{equ-kdv} and estimating the resulting nonlinear commutator
gives the almost conservation law\footnote{As $\sigma\downarrow0$,
the error vanishes and $\|\cdot\|_{G^\sigma}$ tends to the $L^2$ norm, thus recovers
the $L^2$ conservation law.}
\begin{align}\label{equ-into-conserv}
 \sup_{0\leq t\leq\delta}\norm{q(t)}_{G^\sigma}^2
\leq \norm{q_0}_{G^\sigma}^2
+C\sigma^\alpha\norm{q_0}_{G^\sigma}^3,
\qquad 0<\alpha<\tfrac34.	 
\end{align} 
On a bounded ball in $G^\sigma$, the local time step $\delta$ can
be chosen independently of $\sigma$. Using \eqref{equ-into-conserv} repeatedly to reach large time $T$ gives a lower bound   $\sigma_q(T)\geq cT^{-1/\alpha}$,
which improves as $\alpha$ increases. Tesfahun \cite{Tes} used dyadic
space--time bilinear estimates to reach $\alpha=3/4$, giving the bound
$\sigma_q(t)\geq c|t|^{-4/3}$.

In joint work with Huang \cite{HW}, the author obtained the lower
bound $c|t|^{-1/4}$ by introducing a  higher order almost
conservation law in Gevrey spaces (correspongding to \eqref{equ-into-conserv} with $\alpha=4$), which is inspired by the I-method \cite{CKSTT}.
\begin{center}
{\small
\renewcommand{\arraystretch}{1.20}
\begin{tabular*}{0.97\linewidth}{@{\extracolsep{\fill}}lll@{}}
\toprule
\textbf{Work} & \textbf{Analytic radius bound} & \textbf{Method} \\
\midrule
Kato--Masuda \cite{KM} & $\sigma_q(t)\geq\rho_T>0$ for $|t|\leq T$ & Lyapunov functions \\
Bona--Gruji\'c \cite{BG} & $\sigma_q(t)\geq e^{-ct^2}$ & Gevrey energy estimates \\
Bona--Gruji\'c--Kalisch \cite{BGK} & $\sigma_q(t)\geq c|t|^{-12}$ & Bourgain-Gevrey  estimates \\
Selberg--da Silva \cite{SdS} & $\sigma_q(t)\geq c_{\eps}|t|^{-(4/3+\eps)}$ & Almost conservation \\
Tesfahun \cite{Tes} & $\sigma_q(t)\geq c|t|^{-4/3}$ & Sharper energy increment \\
Huang--Wang \cite{HW} & $\sigma_q(t)\geq c|t|^{-1/4}$ & Higher-order corrections \\
\bottomrule
\end{tabular*}
}
\captionsetup{hypcap=false,font=small}
\captionof{table}{Previous lower bounds for the analyticity radius of KdV on
$\R$. }
\label{tab-radius-progress}
\end{center}

In summary, the above earlier quantitative lower bounds tend to zero as time goes to infinity, but they leave open the possibility
that the actual radius stays uniformly positive \cite{HP}. We prove that this is
the case for every analytic datum.

\begin{theorem}\label{thm-fixed-strip}
Let $\sigma_0>0$ and let $q_0\in G^{\sigma_0}(\R)$ be real-valued.
There is a radius
\[
 0<\rho=\rho\bigl(\sigma_0,\norm{q_0}_{G^{\sigma_0}}\bigr)\leq\sigma_0
\]
such that the global solution of \eqref{equ-kdv} satisfies
\begin{equation}\label{equ-intro-10-27}
 q\in C(\R;G^\rho(\R)),\qquad
 \sup_{t\in\R}\norm{q(t)}_{G^\rho}\leq 3\norm{q_0}_{G^{\sigma_0}}.
\end{equation}
\end{theorem}

\begin{remark}\label{rem-radius-loss}One may choose
$$
		\rho\geq c_{\sigma_0}\bigl(1+\norm{q_0}_{G^{\sigma_0}}\bigr)^{-2/3},
$$
	where $c_{\sigma_0}>0$ depends only on $\sigma_0$.
	The radius and norm bound can be chosen uniformly on bounded subsets of
	$G^{\sigma_0}$.
The strip that persists for all time in Theorem~\ref{thm-fixed-strip}
may have to be narrower than the initial strip. For the solutions $q_N$
constructed in Section~\ref{section-multisoliton}, choose
$3\sigma/N<\sigma_0<2\sigma$. Then $q_N(0)\in G^{\sigma_0}(\R)$,
while $q_N(t)\notin G^{\sigma_0}(\R)$ for all sufficiently large $t>0$.
\end{remark}

There are two results on the uniformly positive analytic radius for the KdV equation in different settings.

For the damped KdV equation $u_t+u_{xxx}+uu_x+a(x)u=0$ on $\R$,
with uniformly analytic $a(x)\geq\lambda>0$ and analytic initial
data, Wang \cite{W22} proved that the solution has a fixed
positive analytic radius for all $t\geq0$. In this case, the conservation law \eqref{equ-into-conserv} becomes 
\[
 \norm{u(\delta)}_{G^\sigma}^{2}
 \leq e^{-2\lambda\delta}\norm{u_0}_{G^\sigma}^{2}
   + C\sigma^\alpha\norm{u_0}_{G^\sigma}^{3}
   + \mbox{negligible terms}.
\]
The strict contraction absorbs the terms small in $\sigma$, and
exponential $L^2$ decay makes the remaining errors summable, so a
single sufficiently small $\sigma$ works throughout the iteration.

For the KdV equation on the torus, Kappeler and P\"oschel \cite{KP} proved a uniform fixed positive analytic radius for all time. Their strategy is as follows. If $q(t,x)$ solves the KdV equation, then the spectrum of the Hill operator $-\partial_x^2+q(t,\cdot)$
is independent of time, so its gap lengths of eigenvalues are conserved. For real
periodic potentials, exponential decay of the Fourier coefficients
implies exponential decay of the gap lengths, and the converse holds
with a possibly smaller exponential rate, see e.g. \cite{Pos}. In this way, one obtains a uniform positive analytic width for solutions in all time, see
\cite[Proposition~7]{KP}. However, the method does not work on the real line since a corresponding reconstruction  would involve
the continuous spectral data of the Hill operator.

The proof of Theorem \ref{thm-fixed-strip} goes as follows. We shall use the integrability of the KdV equation in other way. For sufficiently smooth solutions to \eqref{equ-kdv}, we have the following conservation law,
\[
 \cH_m(q)=\frac12\norm{\partial_x^m q}_2^2
            +\sum_{k=3}^{m+2}\cH_{m,k}(q),\qquad
 \cH_m(q(t))=\cH_m(q(0)),\qquad m\geq0,
\]
see e.g.
\cite[Eq.~(3.17) and Corollary~3.8]{KKL}.
Here $\cH_{m,k}$ is the homogeneous term of degree $k$ in $q$ and its
derivatives. Let $h_{m,k}$ be its symmetric Fourier multiplier on
$\Sigma_k=\{\xi_1+\cdots+\xi_k=0\}$ and put
$\mu_k(\xi)=\frac12\sum_{j=1}^k|\xi_j|$. For fixed $k$ and $\xi$,
Theorem~\ref{thm-positive-measure} represents the normalized multipliers
for all $m\geq k-2$ as even moments of a single probability measure.
In particular,
\begin{align}\label{equ-intro-main-esti}
 0\leq h_{m,k}(\xi)
\leq\frac{\Cat_{k-1}}2\binom{2m}{2k-4}
\mu_k(\xi)^{2(m-k+2)},\qquad 3\leq k\leq m+2,
\end{align} 
where $\Cat_{k-1}$ denotes the $(k-1)$st Catalan number, see \eqref{equ-meas-10-30} for a definition.

Set  $a_m(\sigma)=(2\sigma)^{2m}/(2m)!$ and introduce
\[
 \cE_\sigma(q)=\sum_{m\geq0}a_m(\sigma)\cH_m(q),\qquad
 \cQ_\sigma(q)=\frac12\int_\R\cosh(2\sigma\xi)|\wh q(\xi)|^2
                    \frac{\dd\xi}{2\pi}
 \approx\norm q_{G^\sigma}^2.
\]
The simplex structure of the representing measure gives decay in
$k-2$ input frequencies after summation in $m$, as in
\eqref{equ-energy-11-13}. This yields
\[
 |\cE_\sigma(q)-\cQ_\sigma(q)|
 \lesssim \sigma^{3/2}\norm q_{G^\sigma}^3
 \qquad\text{when }\sigma^{3/2}\norm q_{G^\sigma}\text{ is small}.
\]
Choose $0<\rho\leq\sigma_0$ so that $\rho^{3/2}\norm{q_0}_{G^\rho}$
is sufficiently small. It follows that
\[
 \norm{q(t)}_{G^\rho}^2\lesssim\cE_\rho(q(t))
 =\cE_\rho(q_0)\lesssim\norm{q_0}_{G^\rho}^2
\]
throughout the interval of analytic existence. The local theory then
continues the solution at the same width $\rho$. Thus we obtain
$\inf_{t\in\R}\sigma_q(t)\geq\rho>0$.

In Section~\ref{section-positive-measure}, we derive the positive measure
representation for the Hamiltonian multipliers.
With this in hand, in Section~\ref{section-conserved-energy}, we construct the conserved
energy and prove Theorem \ref{thm-fixed-strip}.

\section{Positive measure representation}\label{section-positive-measure}

\subsection{Hamiltonians and their Fourier multipliers}\label{subsection-hamiltonian-multipliers}

In this section, we first recall the classical Riccati construction of the KdV Hamiltonians
\cite{MGK,GD}. In the sequel, let $m$ denote the Hamiltonian
order and $k$ the homogeneous degree.  All series in the auxiliary
parameter $\zeta$ are formal and are used coefficientwise.

For a real Schwartz function $q=q(x)$, let
\begin{equation}\label{equ-meas-10-6}
 P(q;\zeta)=\sum_{j\geq1}p_j(q)\zeta^j
\end{equation}
be the unique differential polynomial series satisfying
\begin{equation}\label{equ-meas-11-29}
 q=\zeta^{-1}P-P_x+P^2.
\end{equation}
Inserting \eqref{equ-meas-10-6} into \eqref{equ-meas-11-29} and comparing the coefficients of $\zeta^j$ on both sides, we have
\begin{equation}\label{equ-meas-11-14}
 p_1=q,\qquad
 p_{j+1}=\partial_xp_j-\sum_{\substack{a+b=j\\a,b\geq1}}p_ap_b,
 \qquad j\geq1.
\end{equation}
For every integer $m\geq0$, define
\begin{equation}\label{equ-meas-11-4}
 \cH_m(q):=\frac{(-1)^{m+1}}2\int_\R p_{2m+3}(q)\dd x.
\end{equation}
In particular, \eqref{equ-meas-11-14} and \eqref{equ-meas-11-4} give 
\begin{align}
 \cH_0(q)&=\frac12\int_\R q^2\dd x,\label{equ-meas-10-19}\\
 \cH_1(q)&=\int_\R\left(\frac12q_x^2+q^3\right)\dd x,
 \label{equ-meas-11-7}
\end{align}
which are the first two conservative quantities for the KdV equation \eqref{equ-kdv}. 

For a polynomial or formal power series $A(s)=\sum_{\ell\geq0}a_\ell s^\ell$,
the notation $[s^k]A(s)$ means the coefficient $a_k$. Thus
$p_{j,k}(q):=[s^k]p_j(sq)$ is the part of $p_j(q)$ homogeneous of
degree $k$ in $q$ and its derivatives. For instance,
$p_3(sq)=s q_{xx}-s^2q^2$, so $p_{3,2}(q)=-q^2$. For every integer $k$, set
\begin{equation}\label{equ-meas-10-17}
 P^{[k]}(q;\zeta):=[s^k]P(sq;\zeta)
                  =\sum_{j\geq1}p_{j,k}(q)\zeta^j.
\end{equation}
Writing
$R_\zeta:=\zeta(1-\zeta\partial_x)^{-1}$, then equation
\eqref{equ-meas-11-29} gives
\begin{equation}\label{equ-meas-10-11}
 P^{[1]}=R_\zeta q,\qquad
 P^{[k]}=-R_\zeta\sum_{a=1}^{k-1}P^{[a]}P^{[k-a]},\quad k\geq2.
\end{equation}
For $3\leq k\leq m+2$, define
\begin{equation}\label{equ-meas-10-23}
 \cH_{m,k}(q):=[s^k]\cH_m(sq)
 =\frac{(-1)^{m+1}}2
       \int_\R[\zeta^{2m+3}]P^{[k]}(q;\zeta)\dd x.
\end{equation}
We write $H^\infty(\R;\R):=\bigcap_{N\geq0}H^N(\R;\R)$ with its
natural Fr\'echet topology.  On this space, the integrated polynomials
are understood by continuous extension after discarding total derivatives.

\begin{proposition}\label{prop-hamiltonian-decomposition}
For every $m\geq0$ and real $q\in H^\infty(\R)$,
\begin{equation}\label{equ-meas-11-15}
 \cH_m(q)=\frac12\norm{\partial_x^mq}_2^2
              +\sum_{k=3}^{m+2}\cH_{m,k}(q),
\end{equation}
where the sum vanishes when $m=0$.  Every monomial in $p_{2m+3,k}$ has
$k$ factors and total differential order $2(m+2-k)$, for
$3\leq k\leq m+2$.
Moreover, if $q\in C(I;H^\infty)$ solves the KdV equation \eqref{equ-kdv} on an open
interval $I$, then
\begin{equation}\label{equ-meas-9-20}
 \frac{\dd}{\dd t}\cH_m(q(t))=0.
\end{equation}
\end{proposition}

\begin{proof}
The result is classical, we only give a sketch for the reader's convenience. By \eqref{equ-meas-11-14}, every monomial in $p_{j,k}$ satisfies
$$
 2k+r=j+1.
$$
Thus $r=2(m+2-k)$ and $k\leq m+2$ when $j=2m+3$.

For the linear term $p_{2m+3,1}=\partial_x^{2m+2}q$, its integral on $\R$ equals to zero. For the quadratic term, since $\int R_\zeta f=\zeta\int f$, using
\eqref{equ-meas-10-11} and integration by parts we get
\[
 \begin{aligned}
 \int_\R p_{2m+3,2}(q)\dd x
 &=-\sum_{a=0}^{2m}\int_\R
              \partial_x^a q\,\partial_x^{2m-a}q\dd x\\
 &=-(-1)^m\sum_{a=0}^{2m}(-1)^a\norm{\partial_x^m q}_2^2
   =(-1)^{m+1}\norm{\partial_x^m q}_2^2.
 \end{aligned}
\]
This proves the decomposition in \eqref{equ-meas-11-15}.

We next show that $\cH_m$ is conserved by the KdV flow. Thanks to 
\cite[Corollary~3.8]{KKL}, the Hamiltonians
\eqref{equ-meas-11-4} Poisson commute with respect to the Gardner Poisson structure.  If
$q\in C(I;H^\infty)$ solves KdV, then $q\in C^1(I;H^{2m+2})$ and
$q_t=\partial_x(\delta\cH_1/\delta q)$. The chain rule gives
$$
 \frac{\dd}{\dd t}\cH_m(q(t))
=\mathrm D\cH_m(q(t))[q_t(t)]
=\int_\R\frac{\delta\cH_m}{\delta q}\,
             \partial_x\frac{\delta\cH_1}{\delta q}\dd x=\{\cH_m,\cH_1\}(q(t))=0.
$$
This proves \eqref{equ-meas-9-20}.
\end{proof}

To sum the hierarchy analytically, we keep track of the Hamiltonian
order $m$ and the homogeneous degree $k$ separately. In other words, we need to establish further quantitative bounds of $\cH_{m,k}(q)$ in terms of $m,k$. To this end, we write each nonlinear term as a Fourier multiplier.

For $k\geq2$, let
\[
 \Sigma_k:=\{\xi\in\R^k:\xi_1+\cdots+\xi_k=0\},
 \qquad
 \dd\Sigma_k:=(2\pi)^{-(k-1)}\dd\xi_1\cdots\dd\xi_{k-1},
 \quad \xi_k=-\sum_{j<k}\xi_j.
\]
For a measurable multiplier $\mathfrak m: \R^k\to \R$ and Schwartz functions
$f_1,\ldots,f_k$, set
\begin{equation}\label{equ-meas-11-6}
 \Lambda_k(\mathfrak m;f_1,\ldots,f_k)
 :=\int_{\Sigma_k}\mathfrak m(\xi)
   \prod_{j=1}^k\wh f_j(\xi_j)\dd\Sigma_k.
\end{equation}

If $\mathfrak m$ is real-valued and satisfies
$\mathfrak m(-\xi)=\mathfrak m(\xi)$, then this integral is real-valued for real-valued $f_1,\cdots,f_k$.
The notation $\Lambda_k(\mathfrak m;q^k)$ denotes $\Lambda_k(\mathfrak m;f_1,\ldots,f_k)$ with $f_1=\cdots=f_k=q$.
Let $S_k$ denote the group of all permutations of $\{1,\ldots,k\}$.
For $\pi\in S_k$, write
$\pi\xi=(\xi_{\pi(1)},\ldots,\xi_{\pi(k)})$.  Define
\begin{equation}\label{equ-meas-9-26}
 [\mathfrak m]_{\Sym}(\xi)
 :=\frac1{k!}\sum_{\pi\in S_k}
 \mathfrak m(\xi_{\pi(1)},\ldots,\xi_{\pi(k)}).
\end{equation}
This normalized symmetrization leaves $\Lambda_k(\mathfrak m;q^k)$ unchanged.
For instance, by the Fourier inversion, we have
$$
 \int_\R\prod_{j=1}^k\partial_x^{r_j}q\dd x
 =\Lambda_k\left(\prod_{j=1}^k(i\xi_j)^{r_j};q^k\right).
$$
Let $h_{m,k}$ be the symmetrized Fourier symbol of
$\frac{(-1)^{m+1}}2p_{2m+3,k}$ restricted to $\Sigma_k$, so that
\begin{equation}\label{equ-meas-10-9}
 \cH_{m,k}(q)=\Lambda_k(h_{m,k};q^k),
 \qquad 3\leq k\leq m+2.
\end{equation}
For $k\geq2$, define the function $\mu_k:\Sigma_k\to[0,\infty)$ by
\begin{equation}\label{equ-meas-9-29}
 \mu_k(\xi):=\frac12\sum_{j=1}^k|\xi_j|,
 \qquad \xi\in\Sigma_k.
\end{equation}
Every nonempty proper subset $I\subset\{1,\ldots,k\}$ satisfies
\begin{equation}\label{equ-meas-10-2}
 \left|\sum_{j\in I}\xi_j\right|
 \leq\min\left\{\sum_{j\in I}|\xi_j|,
                    \sum_{j\notin I}|\xi_j|\right\}
 \leq\mu_k(\xi),  \qquad \xi\in\Sigma_k.
\end{equation}
The first inequality uses $\xi_1+\cdots+\xi_k=0$.  This elementary
localization, rather than an estimate by the full $\ell^1$ mass, is
important for our purpose in this paper.

We use the Catalan numbers, which are defined by induction as
\begin{equation}\label{equ-meas-10-30}
 \Cat_0=1,\qquad
 \Cat_{k-1}=\sum_{a=1}^{k-1}\Cat_{a-1}\Cat_{k-a-1}
 =\frac1k\binom{2k-2}{k-1}\quad(k\geq2).
\end{equation}
The quadratic resolvent recursion generates total mass exactly
$\Cat_{k-1}$ at degree $k$. The key to the proof of Theorem \ref{thm-fixed-strip} is the following result, which is of independent interest.

\begin{theorem}\label{thm-positive-measure}
For every integer $k\geq3$ and every $\xi\in\Sigma_k$, there is a
probability measure $\nu_{k,\xi}$ supported in
$[-\mu_k(\xi),\mu_k(\xi)]$ such that  for every integer
$m\geq k-2$,
\begin{equation}\label{equ-meas-11-27}
  h_{m,k}(\xi)
  =\frac{\Cat_{k-1}}{2}\binom{2m}{2k-4}
    \int_\R y^{2(m-k+2)}\dd\nu_{k,\xi}(y).
\end{equation}
\end{theorem}
\begin{remark}
The measures can be chosen so that
\[
\nu_{k,\pi\xi}=\nu_{k,\xi},
\qquad \pi\in S_k.
\]
For fixed $k$ and $\xi$, this single measure is independent of the
Hamiltonian order $m$.
\end{remark}

Theorem \ref{thm-positive-measure} has the following  direct consequence mentioned in introduction. For $k\geq3$ and
$m\geq k-2$,
\begin{equation}\label{equ-meas-9-17}
 0\leq h_{m,k}(\xi)
 \leq\frac{\Cat_{k-1}}2\binom{2m}{2k-4}
       \mu_k(\xi)^{2(m-k+2)}.
\end{equation}
Here we set $\mu_k^0(\xi)=1$ if $\xi=0$.  In particular,
\begin{equation}\label{equ-meas-9-15}
 h_{k-2,k}(\xi)=\frac{\Cat_{k-1}}2.
\end{equation}

It should be noted that the nonnegativity in \eqref{equ-meas-9-17} only concerns the
multiplier.  For example, $h_{1,3}(\xi)=1$, whereas
$\cH_{1,3}(q)=\int_\R q^3\dd x$ can take either sign.

\subsection{Construction of the representing measure}\label{subsection-measure-construction}

In this subsection, we prove Theorem \ref{thm-positive-measure}. To this end, we first expand the Riccati multipliers into resolvent products and then
express each product as a simplex average.

For $k\geq3$, write $K=\{1,\ldots,k\}$ and define
\[
 S_A(\xi):=\sum_{j\in A}\xi_j,\qquad
 \varnothing\ne A\subseteq K,\quad \xi\in\R^k.
\]
For a nonempty set $A=\{j_1,\ldots,j_a\}\subseteq K$, write
$\xi_A=(\xi_j)_{j\in A}$ and let $M_A(\zeta;\xi_A)$ be the Fourier
symbol of the degree-$a$ term $P^{[a]}(q;\zeta)$, symmetrized in its
$a$ frequency variables $\xi_j$ as in \eqref{equ-meas-9-26}.  Thus for every $\eta\in \R$,
\[
 \widehat{P^{[a]}(q;\zeta)}(\eta)
 =\frac{1}{(2\pi)^{a-1}}
 \int_{\R^{a-1}} M_A(\zeta;\xi_A)
 \prod_{j\in A}\wh q(\xi_j)\,
 \dd\xi_{j_1}\cdots\dd\xi_{j_{a-1}},
 \qquad
 \xi_{j_a}=\eta-\sum_{r=1}^{a-1}\xi_{j_r}.
\]
 Note that when $\xi\in\Sigma_k$ we have $S_K=0$.
But if $A$ is a proper subset of $K$, then we may have $S_A\ne 0$.

Set $z=(i\zeta)^{-1}$ and define $F_A$ by
\[
 M_A(\zeta;\xi_A)=-iF_A\bigl((i\zeta)^{-1}\bigr),
 \qquad\text{equivalently}\qquad
 F_A(z)=iM_A\bigl(1/(iz);\xi_A\bigr).
\]
The factor $-i$ is
chosen to match the Fourier multiplier of $R_\zeta$ in
\eqref{equ-meas-10-11}: at frequency $S_A$ it is
$\zeta/(1-i\zeta S_A)=-i/(z-S_A)$.

\begin{lemma}\label{lem-resolvent-expansion}
Let $\varnothing\ne A\subseteq K$, $a=|A|$  and $\xi\in\Sigma_k$. Then $F_A$ satisfies the following properties.
\begin{description}
\item [(i)] If $a\geq 2$, then the function $(z-S_A(\xi))F_A(z)$ is a finite positive linear
combination of products of exactly $2a-2$ factors
$(z-S_B(\xi))^{-1}$ indexed by nonempty proper subsets
$B\subsetneq A$. If $a=1$, then the factor is $1$.
\item [(ii)] The sum of coefficients in this expansion equals to $\Cat_{a-1}$, namely
\begin{equation}\label{equ-meas-11-26}
	\lim_{z\to\infty}z^{2a-1}F_A(z)=\Cat_{a-1}.
\end{equation}	   
\end{description}
 
\end{lemma}

\begin{proof}
Fix $\xi\in\Sigma_k$. We write $S_B(\xi)$ by $S_B$ for brevity in the proof. 

If $a=1$, from
$P^{[1]}=R_\zeta q$ we obtain $F_{\{j\}}(z)=(z-\xi_j)^{-1}$.

If $a\geq2$, consider the expression $-R_\zeta(P^{[b]}P^{[a-b]})$ in
\eqref{equ-meas-10-11}. For a set $E\subsetneq A$ with $|E|=b$, assign
the labels in $E$ to Fourier transform variables of the first factor and those in $A\setminus E$ to
the second. Since the sum of the Fourier variables is $S_A$, the outer resolvent
$R_\zeta$ acts at frequency $S_A$. It follows that, the contribution
of this assignment is
\[
 -\left(\frac{-i}{z-S_A}\right)(-iF_E(z))(-iF_{A\setminus E}(z))
 =-i\,\frac{F_E(z)F_{A\setminus E}(z)}{z-S_A}.
\]
We also note that for fixed $b$, the $b$ labels in $E$ and the $a-b$ labels in its
complement can each be permuted without changing the respective
symbols. Thus each $E$ occurs $b!(a-b)!$ times among the $a!$
permutations. These considerations and the recursion \eqref{equ-meas-10-11} give
\begin{equation}\label{equ-recur-F}
 F_A(z)=\frac{1}{z-S_A}
\sum_{b=1}^{a-1}\frac{1}{\binom{a}{b}}
\sum_{\substack{E\subsetneq A\\ |E|=b}}
F_E(z)F_{A\setminus E}(z),\qquad a\geq2.
\end{equation}

\medskip 

We prove (i) by induction on $a$. The case $a=1$ follows from the
formula of $F_{\{j\}}$. When $a=2$ and $A=\{i,j\}$, the two assignments give
the same product and hence
\[
 F_{\{i,j\}}(z)=
 \frac{1}{(z-\xi_i-\xi_j)(z-\xi_i)(z-\xi_j)}.
\]
Assume now that (i) holds for every smaller set. Including the outer
resolvents, a term in $F_E F_{A\setminus E}$ has $2b-1$ factors from
nonempty subsets of $E$ and $2(a-b)-1$ factors from nonempty subsets of
$A\setminus E$. Thus we obtain $2a-2$ factors, all indexed by proper
subsets of $A$. Hence, the recursion \eqref{equ-recur-F} combines their products with positive
coefficients, which proves (i).  

\medskip 

To prove (ii), let $c_A$ denote the sum of coefficients  in this expansion.
Clearly, we have $c_{\{j\}}=1$. Note that including an outer resolvent does not change the
coefficients, so the product $F_E F_{A\setminus E}$ also has
coefficient sum $c_Ec_{A\setminus E}$. For $|E|=b$, the induction
hypothesis gives $c_Ec_{A\setminus E}=\Cat_{b-1}\Cat_{a-b-1}$.
With these in hand, averaging over $E$ and summing over $b$, we obtain
\[
 c_A
 =\sum_{b=1}^{a-1}\frac1{\binom{a}{b}}
   \sum_{\substack{E\subsetneq A\\ |E|=b}}c_Ec_{A\setminus E}
 =\sum_{b=1}^{a-1}\Cat_{b-1}\Cat_{a-b-1}
   \overset{\eqref{equ-meas-10-30}}{\xlongequal{\phantom{x}}}\Cat_{a-1}.
\]
Each term of $F_A$ has $2a-1$ resolvent factors including the outer
one. Since $z/(z-S_B)\to1$ for every $B$, the limit of $z^{2a-1}F_A(z)$
is precisely $c_A$. This proves (ii).
\end{proof}

To relate $F_K$ to $h_{m,k}$, observe that integrating in $x$ equals to the Fourier transform at frequency zero, so the input frequencies are on the hyperplane $\Sigma_k$.
Thus, by \eqref{equ-meas-10-23} and the definition of $F_K$, recall that  $z^{-j}=(i\zeta)^j$, we have
\begin{equation}\label{equ-meas-9-11}
 h_{m,k}(\xi)
 =\frac{(-1)^{m+1}}2[\zeta^{2m+3}]
   \bigl[-iF_K((i\zeta)^{-1})\bigr]
 =\frac{1}{2}[z^{-(2m+3)}]F_K(z).
\end{equation}

\begin{lemma}\label{lem-simplex-average}
Let $s_1,\ldots,s_N\in J=[-\mu,\mu]$. Let $\theta$ be a uniformly
distributed random variable vector on the simplex $\theta_i\geq0$, $\sum_{i=1}^N\theta_i=1$,
and put $Y=\sum_{i=1}^N\theta_i s_i$. Then the expectation satisfies
\begin{equation}\label{equ-meas-11-16}
 \mathbb E[(z-Y)^{-N}]=\prod_{i=1}^N(z-s_i)^{-1},
 \qquad z\in\C\setminus J.
\end{equation}
For every $1\leq r<N$, the joint density of any $r$ variables of $\theta_j$ is
\begin{equation}\label{equ-meas-9-25}
 \frac{(N-1)!}{(N-r-1)!}(1-t_1-\cdots-t_r)^{N-r-1}
\end{equation}
on $\{t_i\geq 0:\ t_1+\cdots+t_r\leq 1\}$.
\end{lemma}

\begin{proof}
In case $N=1$, \eqref{equ-meas-11-16} holds clearly. Assume $N\geq2$. Write $\theta'=(\theta_1,\ldots,\theta_{N-1})$ and set
$$ 
 \Delta_d(a):=\{v\in[0,\infty)^d:\textstyle\sum_{i=1}^d v_i\leq a\}, \qquad d\geq0,\ a>0. 
$$
Clearly, the volume $|\Delta_d(a)|=a^d/d!$. Thus the density of $\theta'$ on $\Delta_{N-1}(1)$ is $(N-1)!$.
For real $z>\mu$, the substitution $u=\sum_i u_i$, $u_i=u\theta_i$
gives $\dd u_1\cdots\dd u_N=u^{N-1}\dd u\dd\theta'$. Thus
\[
 \begin{aligned}
  \mathbb E[(z-Y)^{-N}]
  &=(N-1)!\int_{\Delta_{N-1}(1)}(z-Y)^{-N}\dd\theta'\\
  &=\int_{\Delta_{N-1}(1)}\int_0^\infty
       u^{N-1}e^{-u(z-Y)}\dd u\dd\theta'\\
  &=\int_{(0,\infty)^N}e^{-\sum_{i=1}^N u_i(z-s_i)}
       \dd u_1\cdots\dd u_N\\
  &=\prod_{i=1}^N\int_0^\infty e^{-u_i(z-s_i)}\dd u_i
   =\prod_{i=1}^N(z-s_i)^{-1}.
 \end{aligned}
\]
This proves \eqref{equ-meas-11-16} for real $z>\mu$. Moreover, since both sides are holomorphic on $\C\setminus J$, 
the unique continuation property gives \eqref{equ-meas-11-16} for general $z\in \C\setminus J$.

By symmetry, it suffices to compute the density $f_r$ of the first
$r$ coordinates. For $t$ in the interior of $\Delta_r(1)$,
\[
 \begin{aligned}
  f_r(t)
  &=(N-1)!\left|\Delta_{N-r-1}\Bigl(1-\textstyle\sum_{i=1}^r t_i\Bigr)\right|\\
  &=\frac{(N-1)!}{(N-r-1)!}
      \left(1-\sum_{i=1}^r t_i\right)^{N-r-1},
 \end{aligned}
\]
which gives \eqref{equ-meas-9-25}.
\end{proof}

\begin{proof}[Proof of Theorem~\ref{thm-positive-measure}]
Fix $\xi\in\Sigma_k$
and set
\[
 \mu:=\mu_k(\xi),\qquad J:=[-\mu,\mu],
 \qquad E^c:=K\setminus E.
\]

\textbf{Step 1.} We construct a probability measure representing
$z^2F_K(z)$. Since $S_K=0$, the recursion \eqref{equ-recur-F} for $F_K$ gives
\[
 F_K(z)
 =
 \frac1z
 \sum_{b=1}^{k-1}\binom{k}{b}^{-1}
 \sum_{\substack{E\subset K\\|E|=b}}
 F_E(z)F_{E^c}(z).
\]
Since $S_E+S_{E^c}=0$, $z=\frac12\bigl((z-S_E)+(z-S_{E^c})\bigr)$, and hence
\begin{equation}\label{equ-meas-11-17}
 2zF_E(z)F_{E^c}(z)
 =\bigl[(z-S_E)F_E(z)\bigr]F_{E^c}(z)
 +F_E(z)\bigl[(z-S_{E^c})F_{E^c}(z)\bigr].
\end{equation}

By Lemma~\ref{lem-resolvent-expansion}, the first term on the
right-hand side of \eqref{equ-meas-11-17} is a finite positive linear combination of products of $ (2|E|-2)+(2|E^c|-1)=2k-3$ resolvent factors.  The second term is the same.
Every factor is of the form
\[
 (z-S_A(\xi))^{-1},
 \qquad \varnothing\ne A\subsetneq K.
\]
By \eqref{equ-meas-10-2}, we have $|S_A(\xi)|\leq\mu$,
so all poles belong to $J$.

We also record the singleton factors occur in these products.
Each term of $F_A$ contains a factor $(z-\xi_j)^{-1}$ for every $j\in A$.
This holds for singletons. The same conclusion also holds for larger sets because each
recursive split partitions $A$. In \eqref{equ-meas-11-17}, the removed
factor is a singleton factor only if the corresponding set is itself
a singleton. Thus, each product retains singleton factors from at
least $k-1$ distinct input labels.

By Lemma~\ref{lem-simplex-average}, each product in
\eqref{equ-meas-11-17} is the expectation of $(z-Y)^{-(2k-3)}$, where
$Y$ is the convex combination of its poles with uniform simplex
weights. Taking the same positive linear combination of the laws
of $Y$,  we obtain a finite positive measure $\rho_{k,\xi}$ supported in
$J$ such that
\begin{equation}\label{equ-meas-9-19}
 z^2F_K(z)
 =
 \int_J(z-y)^{-(2k-3)}\dd\rho_{k,\xi}(y),
 \qquad z\in\C\setminus J.
\end{equation}
Since $(x/(x-y))^{2k-3}\to 1 (x\to+\infty)$ 
holds uniformly for $y\in J$, using \eqref{equ-meas-11-26}
we obtain
\[
 \begin{aligned}
 \rho_{k,\xi}(J)
 &=
 \lim_{x\to+\infty}
 x^{2k-3}
 \int_J(x-y)^{-(2k-3)}\dd\rho_{k,\xi}(y)\\
 &=
 \lim_{x\to+\infty}x^{2k-1}F_K(x)
 =\Cat_{k-1}.
 \end{aligned}
\]

To make the measure symmetric in $\xi$, apply the construction to
each $\pi\xi$, $\pi\in S_k$, and set
\[
 \nu_{k,\xi}
 :=
 \frac1{k!\Cat_{k-1}}
 \sum_{\pi\in S_k}\rho_{k,\pi\xi}.
\]
All terms in the sum have mass $\Cat_{k-1}$ and are supported in $J$,
since $\mu_k$ is symmetric. Hence $\nu_{k,\xi}$ is a probability
measure on $J$ with
\[
 \nu_{k,\pi\xi}=\nu_{k,\xi},
 \qquad \pi\in S_k,
\]
Since $F_K$ is also symmetric, \eqref{equ-meas-9-19} gives
\begin{equation}\label{equ-meas-9-23}
 z^2F_K(z)
 =
 \Cat_{k-1}
 \int_J(z-y)^{-(2k-3)}\dd\nu_{k,\xi}(y).
\end{equation}
In particular, $\nu_{k,\xi}$ is independent of $m$.

\textbf{Step 2.} We complete the proof. For $|z|>\mu$, we have
\[
 \begin{aligned}
 (z-y)^{-(2k-3)}
 &=
 z^{-(2k-3)}
 \left(1-\frac{y}{z}\right)^{-(2k-3)}\\
 &=
 \sum_{n\geq0}
 \binom{n+2k-4}{2k-4}
 y^n z^{-(n+2k-3)}.
 \end{aligned}
\]
The expansion is absolutely and uniformly convergent for $y\in J$. Thus
for every $m\geq k-2$,
\[
 [z^{-(2m+1)}](z-y)^{-(2k-3)}
 =\binom{2m}{2k-4}y^{2(m-k+2)}.
\]
Substituting this coefficient into \eqref{equ-meas-9-23} and using
\eqref{equ-meas-9-11}, we obtain
\[
 \begin{aligned}
 h_{m,k}(\xi)
 &=
 \frac12[z^{-(2m+3)}]F_K(z) =
 \frac12[z^{-(2m+1)}]\bigl(z^2F_K(z)\bigr)\\
 &=
 \frac{\Cat_{k-1}}2
 \binom{2m}{2k-4}
 \int_J y^{2(m-k+2)}\dd\nu_{k,\xi}(y).
 \end{aligned}
\]
This proves \eqref{equ-meas-11-27}, including the case $\mu=0$,
where $\nu_{k,\xi}=\delta_0$ and $y^0=1$.
\end{proof}

\section{Conserved energy and the fixed analyticity strip}\label{section-conserved-energy}

\subsection{Construction of the conserved energy}\label{subsection-energy-construction}

For every $\sigma>0$, by Taylor expansion we have
\begin{equation}\label{equ-energy-9-8}
 \cosh(2\sigma z)=\sum_{m\geq0}a_m(\sigma)z^{2m},
 \qquad a_m(\sigma)=\frac{(2\sigma)^{2m}}{(2m)!}.
\end{equation}
Introduce the quantity
\begin{equation}\label{equ-energy-10-29}
 \cQ_\sigma(q)=\frac12\int_\R\cosh(2\sigma\xi)
              |\wh q(\xi)|^2\frac{\dd\xi}{2\pi}.
\end{equation}
Thanks to the elementary inequalities $\tfrac12 e^{|x|}\leq\cosh x\leq e^{|x|}$, we have
\begin{equation}\label{equ-energy-9-6}
 \frac14\norm q_{G^\sigma}^2\leq\cQ_\sigma(q)
       \leq\frac12\norm q_{G^\sigma}^2.
\end{equation}
For $k\geq3$ and $\xi\in\Sigma_k$, define the summed multiplier
\begin{equation}\label{equ-energy-11-18}
 \mathcal M_{k,\sigma}(\xi):=
 \sum_{m\geq k-2}a_m(\sigma)h_{m,k}(\xi).
\end{equation}
Set $w_\sigma(u):=\int_0^1 e^{-2\sigma t|u|}\dd t$. Then it is easy to see that
$w_\sigma(u)=w_1(\sigma u)$ and
$$
0\leq w_1(u)\leq\min\{1,(2|u|)^{-1}\}, \quad u\ne 0.
$$
In particular, this implies
$\norm{w_\sigma}_2=\sigma^{-1/2}\norm{w_1}_2<\infty$.

\begin{lemma}\label{lem-multiplier-decay}
For every $\sigma>0$, $k\geq3$ and $\xi\in\Sigma_k$, the series
\eqref{equ-energy-11-18} converges and
\begin{equation}\label{equ-energy-11-13}
 0\leq\mathcal M_{k,\sigma}(\xi)
 \leq\frac{\Cat_{k-1}(2\sigma)^{2k-4}}{2(k-2)!}
       e^{2\sigma\mu_k(\xi)}
       \sum_{\substack{I\subseteq\{1,\ldots,k\}\\ |I|=k-2}}
           \prod_{j\in I}w_\sigma(\xi_j).
\end{equation}
\end{lemma}

\begin{proof}
Fix $\sigma>0$ and $\xi\in\Sigma_k$, and write $r=k-2$ and
$\mu=\mu_k(\xi)$. Multiplying the binomial coefficient in
\eqref{equ-meas-11-27} by $a_m(\sigma)$ gives
$(2\sigma)^{2m}/((2r)!(2m-2r)!)$. Summing the Taylor series of $\cosh$
under the integral by monotone convergence, we obtain
\begin{equation}\label{equ-energy-10-26}
 \mathcal M_{k,\sigma}(\xi)
 =\frac{\Cat_{k-1}(2\sigma)^{2r}}{2(2r)!}
     \int_{-\mu}^{\mu}\cosh(2\sigma y)\dd\nu_{k,\xi}(y).
\end{equation}
By the proof of Theorem~\ref{thm-positive-measure}, we have the decomposition
$\nu_{k,\xi}=\sum_\ell\lambda_\ell\operatorname{Law}(Y_\ell)$,
where $\lambda_\ell\geq0$ and $\sum_\ell\lambda_\ell=1$. In each
component $Y=\sum_{i=1}^{2r+1}\theta_i s_i$, the vector $\theta$ is uniformly
distributed on the simplex, $|s_i|\leq\mu$, and at least $r+1$ poles are
singleton poles $s_{i(j)}=\xi_j$ with distinct input labels.

Fix one component and $\epsilon\in\{-1,1\}$. Then it is easy to see 
\[
 \sum_{\epsilon\xi_j>0}\epsilon\xi_j=\mu
 \quad\Longrightarrow\quad
 \#\{j:\epsilon\xi_j>\mu/2\}\leq1.
\]
Choose a set $I$ of $r$ singleton labels with $\epsilon\xi_j\leq\mu/2$.
Then for $j\in I$, we have $|\xi_j|+\epsilon\xi_j=2\max\{\epsilon\xi_j,0\}\leq\mu$, and thus
\[
 \mu-\epsilon Y
 =\sum_{i=1}^{2r+1}\theta_i(\mu-\epsilon s_i)
 \geq\sum_{j\in I}\theta_{i(j)}(\mu-\epsilon\xi_j)
 \geq\sum_{j\in I}\theta_{i(j)}|\xi_j|.
\]
By \eqref{equ-meas-9-25}, the density of $(\theta_{i(j)})_{j\in I}$
is bounded by $(2r)!/r!$ and supported in $[0,1]^I$. Thus
\[
 \begin{aligned}
  \mathbb E[e^{-2\sigma(\mu-\epsilon Y)}]
  &\leq\frac{(2r)!}{r!}\int_{[0,1]^I}
       \prod_{j\in I}e^{-2\sigma t_j|\xi_j|}\dd t\\
  &=\frac{(2r)!}{r!}\prod_{j\in I}w_\sigma(\xi_j).
 \end{aligned}
\]
It follows that
\[
 \begin{aligned}
  e^{-2\sigma\mu}\int\cosh(2\sigma y)\dd\nu_{k,\xi}(y)
  &=\frac12\sum_{\epsilon=\pm1}\sum_\ell\lambda_\ell
       \mathbb E[e^{-2\sigma(\mu-\epsilon Y_\ell)}]\\
  &\leq\frac{(2r)!}{r!}\sum_{\substack{I\subseteq\{1,\ldots,k\}\\|I|=r}}
       \prod_{j\in I}w_\sigma(\xi_j).
 \end{aligned}
\]
Combining this with \eqref{equ-energy-10-26} proves \eqref{equ-energy-11-13}.
\end{proof}

\begin{proposition}\label{prop-analytic-multilinear}
There is an absolute constant $C$ such that for every $\sigma>0$,
$k\geq3$ and $q_1,\ldots,q_k\in G^\sigma(\R)$,
\begin{equation}\label{equ-energy-10-25}
 \int_{\Sigma_k}\mathcal M_{k,\sigma}(\xi)
      \prod_{j=1}^k|\wh q_j(\xi_j)|\dd\Sigma_k
 \leq\frac{C^k\sigma^{3(k-2)/2}}{(k-2)!}
       \prod_{j=1}^k\norm{q_j}_{G^\sigma}.
\end{equation}
\end{proposition}

\begin{proof}
Set $r=k-2$ and $g_j(u)=e^{\sigma|u|}|\wh q_j(u)|$. For each set
$I\subseteq\{1,\ldots,k\}$ with $|I|=r$, let $F_j=w_\sigma g_j$
if $j\in I$ and $F_j=g_j$ otherwise. By Young's inequality and
Cauchy--Schwarz,
\[
 \begin{aligned}
 & (2\pi)^{k-1}\int_{\Sigma_k}\prod_{j=1}^k F_j(\xi_j)\dd\Sigma_k
  =(F_1*\cdots*F_k)(0)\\
  &\quad\leq\prod_{j\in I}\norm{F_j}_1\prod_{j\notin I}\norm{F_j}_2 \leq\sigma^{-r/2}\norm{w_1}_2^r\prod_{j=1}^k\norm{g_j}_2.
 \end{aligned}
\]
Using $2\mu_k(\xi)=\sum_j|\xi_j|$ and
$\norm{g_j}_2=\sqrt{2\pi}\norm{q_j}_{G^\sigma}$ in
\eqref{equ-energy-11-13}, we obtain
\[
 \begin{aligned}
  \int_{\Sigma_k}\mathcal M_{k,\sigma}(\xi)
       \prod_{j=1}^k|\wh q_j(\xi_j)|\dd\Sigma_k &\leq\frac{\Cat_{k-1}(2\sigma)^{2r}}{2\,r!}\binom{k}{2}
       \left(\frac{\norm{w_1}_2}{\sqrt{2\pi\sigma}}\right)^r
       \prod_{j=1}^k\norm{q_j}_{G^\sigma}\\
  & \leq\frac{C^k\sigma^{3r/2}}{r!}\prod_{j=1}^k\norm{q_j}_{G^\sigma},
 \end{aligned}
\]
where we used $\Cat_{k-1}\leq4^{k-1}$.
\end{proof}

\begin{proposition}\label{prop-conserved-energy}
For every $\sigma>0$, the series
\begin{equation}\label{equ-energy-10-24}
 \cE_\sigma(q):=\sum_{m=0}^\infty a_m(\sigma)\cH_m(q)
\end{equation}
defines a continuous real-valued functional on $G^\sigma$, with the
representation
\begin{equation}\label{equ-energy-9-21}
 \cE_\sigma(q)=\cQ_\sigma(q)
       +\sum_{k\geq3}\Lambda_k(\mathcal M_{k,\sigma};q^k).
\end{equation}
If $q\in C(I;G^\sigma)$  is a real solution of the KdV equation \eqref{equ-kdv} on an
open interval $I$,
\begin{equation}\label{equ-energy-11-1}
 \cE_\sigma(q(t))=\cE_\sigma(q(s)),\qquad s,t\in I.
\end{equation}
\end{proposition}

\begin{proof}
Fix $\sigma>0$. 
Set $A_{k,\sigma}=C^k\sigma^{3(k-2)/2}/(k-2)!$. By \eqref{equ-meas-10-9}
and \eqref{equ-energy-10-25}, we have
\[
 \begin{aligned}
 \sum_{m\geq k-2}a_m(\sigma)|\cH_{m,k}(q)|
 &\leq\int_{\Sigma_k}\mathcal M_{k,\sigma}(\xi)
       \prod_{j=1}^k|\wh q(\xi_j)|\dd\Sigma_k
 \leq A_{k,\sigma}\norm q_{G^\sigma}^k.
 \end{aligned}
\]
The sum $\sum_{k\geq3}$ in \eqref{equ-energy-9-21} converges absolutely and uniformly on
$\norm q_{G^\sigma}\leq R$, since
\[
 \sum_{k\geq3}A_{k,\sigma}R^k
 =C^2R^2\bigl(e^{C\sigma^{3/2}R}-1\bigr).
\]
By \eqref{equ-meas-11-15} and \eqref{equ-energy-9-8}, the quadratic
series sum equals to $\cQ_\sigma(q)$, we have
\[
 \begin{aligned}
 \cE_\sigma(q)
 &=\cQ_\sigma(q)+\sum_{k\geq3}\sum_{m\geq k-2}
       a_m(\sigma)\cH_{m,k}(q)\\
 &=\cQ_\sigma(q)+\sum_{k\geq3}
       \Lambda_k(\mathcal M_{k,\sigma};q^k).
 \end{aligned}
\]
For $\norm q_{G^\sigma},\norm p_{G^\sigma}\leq R$, by Cauchy-Schwarz's inequality we have
\[
 \begin{aligned}
 |\cQ_\sigma(q)-\cQ_\sigma(p)|
 &\leq R\norm{q-p}_{G^\sigma},\\
 |\Lambda_k(\mathcal M_{k,\sigma};q^k)
   -\Lambda_k(\mathcal M_{k,\sigma};p^k)|
 &\leq k A_{k,\sigma}R^{k-1}\norm{q-p}_{G^\sigma}.
 \end{aligned}
\]
Hence
\[
 |\cE_\sigma(q)-\cE_\sigma(p)|
 \leq\left(R+\sum_{k\geq3}k A_{k,\sigma}R^{k-1}\right)
      \norm{q-p}_{G^\sigma}.
\]
This proves the continuous of the functional since the sum is finite.

  Let $q\in C(I;G^\sigma)$. Then $q\in C(I;H^\infty)$, and each $\cH_m(q(t))$ is constant as time $t$ varies by
Proposition~\ref{prop-hamiltonian-decomposition}. It follows that for $s,t\in I$,
\[
 \begin{aligned}
 \cE_\sigma(q(t))
 &=\lim_{M\to\infty}\sum_{m=0}^M a_m(\sigma)\cH_m(q(t))\\
 &=\lim_{M\to\infty}\sum_{m=0}^M a_m(\sigma)\cH_m(q(s))
  =\cE_\sigma(q(s)).
 \end{aligned}
\]
This proves \eqref{equ-energy-11-1}.
\end{proof}

\subsection{Coercivity and global continuation}\label{subsection-global-continuation}

As mentioned in the introduction, the local persistence of the KdV equation \eqref{equ-kdv} in $G^\sigma$ follows directly from the local well posedness in $G^\sigma$. The latter can be proved by the
bilinear estimate of Kenig--Ponce--Vega \cite{KPV}, see e.g. \cite{HW,SdS}. More precisely, we have 

\begin{proposition}\label{prop-local-analytic-persistence}
For every $R>0$, there is $\delta_R\in(0,1]$
independent of $\sigma>0$ and $t_0\in\R$ such that for every $\norm{q(t_0)}_{G^\sigma}\leq R$,  \eqref{equ-kdv} has a unique solution $q\in C([t_0-\delta_R,t_0+\delta_R];G^\sigma)$.
\end{proposition}

Based on the conserved energy constructed in Proposition \ref{prop-conserved-energy}, now we show that the solution has a uniform positive analytic radius for all time. 

The representation \eqref{equ-energy-9-21} and
Proposition~\ref{prop-analytic-multilinear} give
\begin{equation}\label{equ-energy-10-10}
 |\cE_\sigma(q)-\cQ_\sigma(q)|
 \leq\norm q_{G^\sigma}^2
       \sum_{k\geq3}\frac{C^k}{(k-2)!}
        \bigl(\sigma^{3/2}\norm q_{G^\sigma}\bigr)^{k-2}.
\end{equation}
This estimate, together with \eqref{equ-energy-9-6}, 
implies $\cE_\sigma(q)\geq 0$ when $\sigma^{3/2}\norm q_{G^\sigma}$ is small.

\begin{proposition}\label{prop-small-data-strip}
There is an absolute $\eta_*>0$ such that, for every $\sigma>0$ and
$q_0\in G^\sigma(\R)$ with $\sigma^{3/2}\norm{q_0}_{G^\sigma}\leq\eta_*$,
the solution constructed in Proposition \ref{prop-local-analytic-persistence} satisfies
\begin{equation}\label{equ-energy-11-11}
 q\in C(\R;G^\sigma),\qquad
 \sup_{t\in\R}\norm{q(t)}_{G^\sigma}
 \le3\norm{q_0}_{G^\sigma}.
\end{equation}
\end{proposition}

\begin{proof}
We split the proof into two steps.
	
\textbf{Step 1.} We bound the analytic norm on its maximal interval
of continuity. First note that we can choose $\eta_*>0$ small enough such that
\[
 \sum_{k\ge3}\frac{C^k}{(k-2)!}(3\eta_*)^{k-2}\le\frac18,
\]
since the series tends to zero as $\eta_*\downarrow 0$.
If $\sigma^{3/2}\norm q_{G^\sigma}\leq3\eta_*$, then
\eqref{equ-energy-10-10} and
\eqref{equ-energy-9-6} give
\begin{align}\label{equ-energy-compar}
 \frac18\norm{q}_{G^\sigma}^2
\le\cE_\sigma(q)\le\frac58\norm{q}_{G^\sigma}^2.
\end{align}

Let $q$ be the solution given in Proposition \ref{prop-local-analytic-persistence} and set $a=\norm{q_0}_{G^\sigma}$.
If $a=0$, the result follows clearly by uniqueness. For $a>0$, let $I_\sigma$
be the maximal open interval containing $0$ on which
$q\in C(I_\sigma;G^\sigma)$. Since $\sigma^{3/2}a\leq\eta_*$,
the comparison \eqref{equ-energy-compar} applies at $t=0$ and at any time $t_*$ with
$\norm{q(t_*)}_{G^\sigma}=3a$. At such a time,
\eqref{equ-energy-11-1} would give
\[
 \frac98a^2
 \le\cE_\sigma(q(t_*))
 =\cE_\sigma(q_0)\le\frac58a^2.
\]
This is impossible. By continuity,
\[
 \sup_{t\in I_\sigma}\norm{q(t)}_{G^\sigma}\leq3a.
\]

\textbf{Step 2.} We prove that $I_\sigma=\R$. Let $\delta_{3a}$ be the
lifespan in Proposition~\ref{prop-local-analytic-persistence}.
If $T:=\sup I_\sigma<\infty$, choose $t_0\in I_\sigma$ such that
\[
 T-\frac{\delta_{3a}}2<t_0<T.
\]
Since $\norm{q(t_0)}_{G^\sigma}\leq3a$, Proposition~\ref{prop-local-analytic-persistence}
gives
\[
 q\in C([t_0-\delta_{3a},t_0+\delta_{3a}];G^\sigma).
\]
Since $t_0+\delta_{3a}>T$, this leads a  contradiction with the maximality of $I_\sigma$. Thus $T=+\infty$.
The argument at the left endpoint is the same. Thus $I_\sigma=\R$,
and the estimate from Step 1 gives \eqref{equ-energy-11-11}.
\end{proof}

Now we apply Proposition~\ref{prop-small-data-strip} to prove Theorem~\ref{thm-fixed-strip}.

\begin{proof}[Proof of Theorem~\ref{thm-fixed-strip}]
Let $N=\norm{q_0}_{G^{\sigma_0}}$ and choose
\begin{equation}\label{equ-energy-11-10}
 \rho:=\min\left\{\sigma_0,\left(\frac{\eta_*}{1+N}\right)^{2/3}\right\}.
\end{equation}
Since $\rho\leq\sigma_0$, we have $\norm{q_0}_{G^\rho}\leq N$ and hence
\[
 \rho^{3/2}\norm{q_0}_{G^\rho}
 \leq \eta_*\frac{N}{1+N}\leq\eta_*.
\]
Applying Proposition~\ref{prop-small-data-strip} with $\sigma=\rho$, we obtain 
$q\in C(\R;G^\rho(\R))$ and
$$
\sup_{t\in\R}\norm{q(t)}_{G^\rho}\leq3\norm{q_0}_{G^\rho}\leq3N.
$$
Moreover, \eqref{equ-energy-11-10} gives the
 lower bound 
$$
	\rho\geq c_{\sigma_0}\bigl(1+\norm{q_0}_{G^{\sigma_0}}\bigr)^{-2/3},
$$ 
 with $c_{\sigma_0}=\min\{\sigma_0,\eta_*^{2/3}\}$.
This completes the proof.
\end{proof}

\section{A multisoliton obstruction}\label{section-multisoliton}

In this section, we show that the infimum of the analyticity radius over time can be
arbitrarily small among solutions with the same initial radius. Thus Theorem \ref{thm-fixed-strip} can not be expected to be improved to $\sigma_q(t)= \sigma_0$.
For a global solution $q$ of \eqref{equ-kdv}, we write
\[
 \rho(q):=\inf_{t\in\R}\sigma_q(t).
\]

\begin{proposition}\label{prop-multisoliton-radius}
Let $\sigma>0$ and let $N\geq2$ be an integer.  Set
$a=\pi/(4\sigma)$. If
\begin{equation}\label{equ-soliton-9-4}
 q_N(0,x)=-a^2N(N+1)\operatorname{sech}^2(ax),
\end{equation}
then the solution $q_N$ of \eqref{equ-kdv} is an $N$-soliton solution and satisfies
\begin{equation}\label{equ-soliton-10-5}
 \sigma_{q_N}(0)=2\sigma,
 \qquad
 \rho(q_N)<\frac{3\sigma}{N}.
\end{equation}
\end{proposition}

Before giving the proof, we make some comments on Proposition \ref{prop-multisoliton-radius}. The profiles in \eqref{equ-soliton-9-4} are classical reflectionless
P\"oschl--Teller potentials. In  \cite[Theorem~2]{BW}, Bona--Weissler
obtained large-time pole asymptotics for two-soliton solutions,
corresponding to the case $N=2$ above. Moreover, McCue et al.~\cite[Section~4]{McCue}
used the same family 
with detailed analyses for $N=2$ and $N=3$. In the sequel, our argument locates
a pole near the fastest soliton at large times for every $N\geq2$,
giving the $O(N^{-1})$ upper bound on $\rho(q_N)$ at fixed initial radius.

For $j\geq0$, define
\begin{equation}\label{equ-soliton-10-12}
 L_j:=-\partial_x^2-a^2j(j+1)\operatorname{sech}^2(ax),
 \qquad
 A_j:=\partial_x+ja\tanh(ax) \quad (j\geq1).
\end{equation}

\begin{lemma}\label{lem-poschl-teller-data}
The negative eigenvalues of $L_N$ are precisely
\begin{equation}\label{equ-soliton-11-28}
 -\kappa_j^2=-a^2j^2,
 \qquad 1\leq j\leq N,
\end{equation}
and $L_N$ is reflectionless.
\end{lemma}

\begin{proof}
We first determine the negative eigenvalues of $L_N$. The formal adjoint
of $A_j$ is $A_j^*=-\partial_x+ja\tanh(ax)$, and direct multiplication
gives
\begin{equation}\label{equ-soliton-9-27}
 A_j^*A_j=L_j+j^2a^2,
 \qquad
 A_jA_j^*=L_{j-1}+j^2a^2,
\end{equation}
and hence
\begin{equation}\label{equ-soliton-9-3}
 A_jL_j=L_{j-1}A_j,
 \qquad
 A_j^*L_{j-1}=L_jA_j^*.
\end{equation}

Since $A_j\operatorname{sech}^j(ax)=0$,
\eqref{equ-soliton-9-27} gives
\[
 L_j\operatorname{sech}^j(ax)=-j^2a^2\operatorname{sech}^j(ax).
\]
Let $m>j$ and suppose $L_{m-1}\phi=-j^2a^2\phi$ with $\phi\ne0$.
By   \eqref{equ-soliton-9-3}, $A_m^*\phi$ solves
the eigenvalue equation for $L_m$ with the same eigenvalue. It is
nonzero since
\[
 \norm{A_m^*\phi}_2^2
 =\langle\phi,A_mA_m^*\phi\rangle
 =(m^2-j^2)a^2\norm{\phi}_2^2>0.
\]
Iterating from $m=j+1$ to $m=N$, we see that
$A_N^*\cdots A_{j+1}^*\operatorname{sech}^j(ax)$ is a nonzero $L^2$
eigenfunction of $L_N$ with eigenvalue $-a^2j^2$.

To rule out other negative eigenvalues, let $L_j\phi=\lambda\phi$
with $\lambda<0$ and $\phi\in L^2(\R)$. By the standard elliptic regularity, we have $\phi\in H^2(\R)$, so $A_j\phi\in L^2(\R)$. If $A_j\phi=0$,
then $\phi$ is a multiple of $\operatorname{sech}^j(ax)$ and
$\lambda=-j^2a^2$. If $A_j\phi\ne0$, the first identity in
\eqref{equ-soliton-9-3} gives an eigenfunction of $L_{j-1}$ with
eigenvalue $\lambda$. Since $L_0=-\partial_x^2$ has no negative
$L^2$ eigenvalues, this procedure must stop in $\ker A_m$ for some
$1\leq m\leq j$. Thus $\lambda=-m^2a^2$, which proves
\eqref{equ-soliton-11-28}.

We next show that $L_N$ is reflectionless.
For real $k\ne0$, set
\[
 \psi_N(k,x):=A_N^*A_{N-1}^*\cdots A_1^*e^{ikx}.
\]
 Using \eqref{equ-soliton-9-3} repeatedly we get 
$L_N\psi_N=k^2\psi_N$.  Since
$$
\tanh(ax)=1+O(e^{-2ax}) \quad (x\to+\infty)
\qquad \tanh(ax)=-1+O(e^{2ax}) \quad (x\to-\infty),
$$
 we have
\begin{equation}\label{equ-soliton-10-14}
 \psi_N(k,x)=
 \begin{cases}
  \displaystyle
  \left(\prod_{j=1}^N(ja-ik)\right)e^{ikx}+o(1),&x\to+\infty,\\[14pt]
  \displaystyle
  \left(\prod_{j=1}^N(-ja-ik)\right)e^{ikx}+o(1),&x\to-\infty.
 \end{cases}
\end{equation}
Normalizing
$\psi_N$ at $+\infty$, we obtain the right Jost solution. Its expansion at
$-\infty$ has zero coefficient of $e^{-ikx}$, so the reflection
coefficient vanishes.
\end{proof}

\medskip
In the sequel, we set
\begin{equation}\label{equ-soliton-10-8}
 \kappa_j=aj,
 \qquad 1\leq j\leq N.
\end{equation}

\begin{lemma}\label{lem-multisoliton-poles}
Let $q_N(0)$ be the potential in \eqref{equ-soliton-9-4}, and let
$\kappa_j$ be given by \eqref{equ-soliton-10-8}.  There are positive constants
$\gamma_1,\ldots,\gamma_N$ such that the  solution $q_N$ of \eqref{equ-kdv} can be
represented by
\begin{equation}\label{equ-soliton-11-23}
 q_N(t,x)=-2\partial_x^2\log\tau_N(x,t),
 \qquad
 \tau_N(x,t)=\det(I+M(x,t)),
\end{equation}
where
\begin{equation}\label{equ-soliton-10-21}
 M_{jk}(x,t)
 =\frac{\gamma_j\gamma_k}{\kappa_j+\kappa_k}
   e^{-(\kappa_j+\kappa_k)x+4(\kappa_j^3+\kappa_k^3)t}.
\end{equation}
For each fixed $t$, the function $\tau_N(\cdot,t)$ is entire and is positive
on the real axis.  If $z_0\notin\R$ is a simple zero of
$\tau_N(\cdot,t)$, then the meromorphic continuation of $q_N(t,\cdot)$ has
a nonremovable double pole at $z_0$, with the local expansion
\eqref{equ-soliton-9-28} below.  In particular,
$\sigma_{q_N}(t)\leq |\operatorname{Im}z_0|$.
\end{lemma}

\begin{proof}
The potential \eqref{equ-soliton-9-4} is real and exponentially decaying.
By Lemma~\ref{lem-poschl-teller-data}, it is reflectionless with bound-state
parameters $\kappa_j=aj$, $1\leq j\leq N$. The finite reflectionless
reconstruction theorem gives  the representation
\eqref{equ-soliton-11-23}--\eqref{equ-soliton-10-21}, see
\cite[formulae~(1.1) and~(1.3)]{HMM}. This
representation solves \eqref{equ-kdv} and recovers
\eqref{equ-soliton-9-4} at $t=0$.

For $x,t\in \R$, set
$v_j(s,t)=\gamma_j e^{-\kappa_js+4\kappa_j^3t}$.  Then by \eqref{equ-soliton-10-21} one can check that
\begin{equation}\label{equ-soliton-11-12}
 M_{jk}(x,t)=\int_x^\infty v_j(s,t)v_k(s,t)\,\dd s.
\end{equation}
Thanks to \eqref{equ-soliton-11-12},  we have
\[
 c^*M(x,t)c
 =\int_x^\infty\left|\sum_{j=1}^N c_jv_j(s,t)\right|^2\dd s
 \geq0, \qquad c\in \C^N.
\]
It follows that $M(x,t)$ is positive semidefinite, and
\[
 \tau_N(x,t)=\det(I+M(x,t))
 =\prod_{j=1}^N\bigl(1+\lambda_j(M(x,t))\bigr)\geq1.
\]
By \eqref{equ-soliton-10-21}, each entry of $M(z,t)$ is entire in $z$,
so $\tau_N(\cdot,t)$ is entire as well.

On the real axis, \eqref{equ-soliton-11-23} may be written as
\[
 q_N(t,x)=-2\partial_x\left(
 \frac{\partial_x\tau_N(x,t)}{\tau_N(x,t)}\right).
\]
The logarithmic derivative $\partial_z\tau_N/\tau_N$ is meromorphic
in $z$. Thus 
\begin{equation}\label{equ-soliton-11-24}
 Q_t(z):=-2\left(\frac{\partial_z\tau_N(z,t)}{\tau_N(z,t)}\right)'
 =-2\frac{\tau_N''(z,t)\tau_N(z,t)-(\tau_N'(z,t))^2}
          {\tau_N(z,t)^2}.
\end{equation}
Let $z_0$ be a simple zero of $\tau_N(\cdot,t)$.  Then near $z_0$, we can write 
\[
 \tau_N(z,t)=(z-z_0)a_t(z),
 \qquad a_t(z_0)\ne0,
\]
with $a_t$ holomorphic near $z_0$.  Substitution this into
\eqref{equ-soliton-11-24} gives
\begin{equation}\label{equ-soliton-9-28}
 Q_t(z)=\frac{2}{(z-z_0)^2}+H_t(z),
\end{equation}
where $H_t$ is holomorphic near $z_0$.  Hence the pole is nonremovable.

If $q_N(t)\in G^s(\R)$ for some $s>|\operatorname{Im}z_0|$, then the
Paley-Wiener theorem gives a holomorphic extension $h$ of
$q_N(t,\cdot)$ to $|\operatorname{Im}z|<s$. Thus the function
\[
 h(z)\tau_N(z,t)^2
 +2\bigl(\tau_N''(z,t)\tau_N(z,t)-(\tau_N'(z,t))^2\bigr)
\]
is holomorphic in this strip and vanishes on the real axis by
\eqref{equ-soliton-11-24}, so it vanishes throughout the strip by the
identity theorem. In particular, letting $z=z_0$ we obtain
$-2(\tau_N'(z_0,t))^2=0$, which leads to a contradiction with simplicity of the zero.
This proves $\sigma_{q_N}(t)\leq |\operatorname{Im}z_0|$.
\end{proof}

We now prove Proposition~\ref{prop-multisoliton-radius}.

\begin{proof}[Proof of Proposition~\ref{prop-multisoliton-radius}]
By Lemma~\ref{lem-multisoliton-poles}, the desired solution is constructed by the $N$-soliton representation. It remains to show the analytic radius estimates.
We first compute it for the initial data. In fact, since $a=\pi/(4\sigma)$ and
$$
 \int_\R \operatorname{sech}^2(ax)e^{-ix\xi}\,\dd x
 =\frac{\pi\xi}{a^2\sinh(\pi\xi/(2a))}
$$
we have
$$
 \widehat{q_N(0,\xi)} 
 =-\frac{\pi N(N+1)\xi}{\sinh(2\sigma\xi)}.
$$
Thus $\widehat{q_N(0)}$ is bounded near $\xi=0$, and
$$
 |\widehat{q_N(0,\xi)}|\approx_{\sigma,N}
 |\xi|e^{-2\sigma|\xi|},
 \qquad |\xi|\geq1.
$$
It follows that $e^{s|\xi|}\widehat{q_N(0,\xi)}\in L^2_\xi(\R)$ if and only if $s<2\sigma$. Therefore, $\sigma_{q_N}(0)=2\sigma$.

\medskip

To bound the analytic radius at later times, we set
\begin{equation}\label{equ-soliton-9-2}
 z=4\kappa_N^2t+y,
 \qquad
 F_t(y):=\tau_N(4\kappa_N^2t+y,t).
\end{equation}
The translation is real, so $\operatorname{Im}z=\operatorname{Im}y$.
Substitution into \eqref{equ-soliton-10-21} gives
\begin{equation}\label{equ-soliton-9-13}
 M_{jk}(4\kappa_N^2t+y,t)
 =\frac{\gamma_j\gamma_k}{\kappa_j+\kappa_k}
  e^{-(\kappa_j+\kappa_k)y}e^{-4t(d_j+d_k)},
 \qquad
 d_j:=\kappa_j(\kappa_N^2-\kappa_j^2).
\end{equation}
Since $d_N=0$ and $d_j>0$ for $j<N$, all entries except the
$(N,N)$ entry tend to zero locally uniformly as $t\to+\infty$. The
remaining entry is
\[
 M_{NN}(4\kappa_N^2t+y,t)
 =\frac{\gamma_N^2}{2\kappa_N}e^{-2\kappa_Ny}.
\]
Taking the determinant yields, locally uniformly in $\C$,
\begin{equation}\label{equ-soliton-9-14}
 F_t(y)\longrightarrow F(y):=1+b_Ne^{-2\kappa_Ny},
 \qquad
 b_N:=\frac{\gamma_N^2}{2\kappa_N}>0,
 \quad t\to+\infty.
\end{equation}

The zeros of $F$ are
\[
 y_\ell=\frac{\log b_N+i(2\ell+1)\pi}{2\kappa_N},
 \qquad \ell\in\mathbb Z.
\]
Let
\begin{equation}\label{equ-soliton-10-13}
 y_*:=\frac{\log b_N+i\pi}{2\kappa_N}.
\end{equation}
Since $\kappa_N=aN=\pi N/(4\sigma)$,
\[
 \operatorname{Im}y_*=\frac{2\sigma}{N},
 \qquad
 |y_\ell-y_*|\geq\frac{4\sigma}{N}\quad(\ell\ne0).
\]
Set
\[
 r_N:=\frac{\sigma}{N},
 \qquad
 D_N:=\{y\in\C:|y-y_*|<r_N\}.
\]
Then $\overline{D_N}$ lies in the upper half-plane and contains no zero of
$F$ other than the simple zero $y_*$.  Hence
\[
 m_N:=\min_{|y-y_*|=r_N}|F(y)|>0.
\]
By the locally uniform convergence in \eqref{equ-soliton-9-14}, there exists
$T_N=T_N(\sigma,N)>0$ such that, for every $t\geq T_N$,
\begin{equation}\label{equ-soliton-9-30}
 \sup_{|y-y_*|=r_N}|F_t(y)-F(y)|<m_N\leq |F(y)|.
\end{equation}
By Rouch\'e's theorem, $F_t$ has exactly one zero $y(t)$ in $D_N$,
counted with multiplicity. Thus $y(t)$ is simple. Set
\[
 z(t):=4\kappa_N^2t+y(t).
\]
Then $\tau_N(z(t),t)=0$, and
\[
 0<\operatorname{Im}z(t)
 =\operatorname{Im}y(t)
 <\operatorname{Im}y_*+r_N
 =\frac{3\sigma}{N},
 \qquad t\geq T_N.
\]
Applying Lemma~\ref{lem-multisoliton-poles} gives
\[
 \sigma_{q_N}(t)<\frac{3\sigma}{N},
 \qquad t\geq T_N.
\]
Consequently,
\begin{equation}\label{equ-soliton-10-3}
 \rho(q_N)<\frac{3\sigma}{N}.
\end{equation}
This proves \eqref{equ-soliton-10-5}.
\end{proof}

\section*{Acknowledgements}
\begin{sloppypar}
This work was partially supported by the National Natural Science Foundation of
China (No.~12571260), the Natural Science Foundation of Hunan Province
(No.~2026JJ20012), and the Hunan Basic Science Research Center for
Mathematical Analysis (2024JC2002).
\end{sloppypar}

\section*{History and AI disclosure}
Shortly after completing the 2019 paper with Huang \cite{HW}, the first
author realized that extending its method to higher-order modified
energies would be impractical: the proof already required extensive
factorization of high-degree polynomials, and each additional
correction would make the algebra substantially more involved.
Inspired by the theory of global attractors for damped KdV, the first
author then found that damping supplied a contraction that could
absorb the errors that vanish as the strip narrows, yielding a fixed
radius on the real line \cite{W22}. A referee of that paper drew the
first author's attention to a result then new to him: Kappeler and
P\"oschel \cite{KP} had obtained a fixed strip for undamped periodic
KdV. The referee suggested that the analogous result should hold on
the real line. We began studying their method systematically in
September 2025 and soon saw that their reconstruction estimates
depended on periodic spectral data and could not be applied
directly on the line. Their use of integrability suggested a way
forward: to obtain quantitative estimates across the hierarchy of
conserved Hamiltonians. In discussions with ChatGPT, we initially
sought a pointwise multiplier bound of the form \eqref{equ-intro-main-esti}
to improve the known lower bounds on the analyticity radius. ChatGPT
ultimately provided the stronger positive measure representation in
Theorem~\ref{thm-positive-measure}. This representation allowed us to
sum the Hamiltonians with analytic weights into a conserved energy
and prove a fixed positive lower bound for the analyticity radius. AI also helped
us develop the counterexample in Section~\ref{section-multisoliton}
and polish the manuscript.

\begingroup
\small
\linespread{0.93}\selectfont
\setlength{\labelsep}{0.4em}

\endgroup

\medskip
\noindent\begin{minipage}{\linewidth}
\small
\textsc{Ming Wang and Jie Zeng}\\
School of Mathematics and Statistics, HNP-LAMA, Central South University\\
Changsha, Hunan 410083, P.R.\ China\\
\textit{Email:} \texttt{m.wang@csu.edu.cn} \qquad \texttt{zengjie@csu.edu.cn} 
\end{minipage}

\end{document}